\documentclass[11pt]{article}
\usepackage{amsfonts}
\usepackage{subfigure} 
\usepackage{epsfig}
\usepackage{amsmath, amsthm}
\newcommand{\be}{\begin{eqnarray}}     	\newcommand{\ee}{\end{eqnarray}}

\newcommand{\vol}{\mathrm{Vol}}

\newcommand{\rem}{\mathrm{Rm}}
\newcommand{\ric}{\mathrm{Ric}}

\newcommand{\Id}{\mathrm{Id}}
 
\title{Convergence of the Ricci flow toward a unique soliton}
\author{Natasa Sesum}
\date{}
\theoremstyle{plain}
\newtheorem{dummy}{Dummy}

\theoremstyle{definition}
\newtheorem{definition}[dummy]{Definition}
\newtheorem{step}{Step}[dummy]
\newtheorem{case}{Case}

\theoremstyle{plain}
\newtheorem{corollary}[dummy]{Corollary}
\newtheorem{remark}[dummy]{Remark}
\newtheorem{lemma}[dummy]{Lemma}
\newtheorem{theorem}[dummy]{Theorem}
\newtheorem{proposition}[dummy]{Proposition}
	
\newtheorem{claim}[dummy]{Claim}

\begin{document}

\maketitle

\begin{abstract}
We will consider a {\it $\tau$-flow}, given by the equation
$\frac{d}{dt}g_{ij} = -2R_{ij} + \frac{1}{\tau}g_{ij}$ on a closed
manifold $M$, for all times $t\in [0,\infty)$. We will prove that if
the curvature operator and the diameter of $(M,g(t))$ are uniformly
bounded along the flow and if one of the limit solitons is integrable,
then we have a convergence of the flow toward a unique soliton, up to
a diffeomorphism.
\end{abstract}

\begin{section}{Introduction}

The Ricci flow equation 
$$\frac{d}{dt}g_{ij} = -2R_{ij},$$ 
has been introduced by R. Hamilton
in his seminal paper \cite{hamilton1982}. We will refer to this
equation as to an unnormalized Ricci flow. A normalized Ricci flow is
given by the equation
$$\frac{d}{dt}\tilde{g}_{ij} = -2R(\tilde{g})_{ij} +
\frac{2}{n}r\tilde{g}_{ij},$$ 
where $r = \frac{1}{\vol(M)}\int_M R(\tilde{g})dV_{\tilde{g}}$. 
This equation is sometimes more
convenient to consider, since a volume of a manifold is being fixed
along the normalized Ricci flow and a volume collapsing case can not
happen in a limit, if the limit exists.

A natural question that arises in studying the evolution equations, in
particular the Ricci flow equation, is under which conditions a
solution will exist for all times, that is under which conditions it
will avoid the singularities at finite times. The other question one
can ask is if there exists a limit of the solutions when we approach
infinity and how we can describe the metrics obtained in the limit. In
the case of dimension three with positive Ricci curvature and
dimension four with positive curvature operator we know (due to
R. Hamilton) that the solutions of the Ricci flow equation, in both
cases exist for all times, converging to Einstein metrics. In general,
we can not expect to get an Einstein metric in the limit. We can
expect to get in the limit a solution to the Ricci flow equation which
moves under one-parameter subgroup of the symmetry group of the
equation. These kinds of solutions are called solitons. Since the
Ricci flow equation is a gradient flow of Perelman's functional
$\mathcal{W}$, it is natural to expect that a soliton in the limit is
unique up to diffeomorphisms.

Our goal in this paper is to prove the following theorem.

\begin{theorem}
\label{theorem-theorem_uniqueness}
Let $(g_{ij})_t = -2R_{ij} + \frac{1}{\tau}g_{ij}$ be a Ricci flow on
a closed manifold $M$ with uniformly bounded curvature operators and
diameters for all $t\in [0,\infty)$. Assume also that some limit
soliton is integrable. Then there is an $1$-parameter family of
diffeomorphisms $\phi(t)$, a unique soliton $h(t)$ and constants $C$,
$\delta$, $t_0$ such that $|\phi(t)^*g(t) - h(0)|_{k,\alpha} <
Ce^{-\delta t}$, for all $t\in [t_0,\infty)$. Moreover, if $\psi(t)$ is
a diffeomorphism such that $h(t) = \psi^* h(0)$, then
$|(\phi\psi)^*g(t) - h(t)|_{C^0} < Ce^{-ct}$.
\end{theorem}

The ideas for the proof of Theorem \ref{theorem-theorem_uniqueness}
have been inspired by those of Cheeger and Tian in \cite{cheeger1994}.
\\
\\
{\it Outline of the proof of Theorem \ref{theorem-theorem_uniqueness}}
\\
\\
In order to deal with this problem, we will first construct a gauge on
time intervals of an arbitrary length, so that in the chosen gauge the
$\tau$-flow equation becomes strongly parabolic. We will look at the
solutions of a strictly parabolic equation. It will turn out that our
metrics (in the right gauge) will satisfy a strictly parabolic
equation that is almost linear and therefore their behavior is
modeled on the behavior of the solutions of the linear
equation. There are $3$ types of the solutions of our strictly
parabolic equation,
\begin{itemize}
\item
the solutions that have an exponential growth,
\item
the solutions that have an exponential decay,
\item
the solutions that change very slowly.  
\end{itemize}
Roughly speaking, the integrability condition means that the solutions
of a linearized deformation equation for solitons arise from a curve
of metrics satisfying the same soliton equation. To deal with those
slowly changing solutions we will use the integrability condition to
change the reference soliton metric so that at the end we deal only
with the cases of either a growth or a decay.  We will rule out the
possibility of the exponential growth, by using the fact that our flow
sequentially converges toward solitons and by using the similar
arguments established by L.Simon in \cite{simon1983} and also later
used by Cheeger and Tian in \cite{cheeger1994}. We will be left with
the exponential decay which will allow us to continue our gauge up to
infinity.

The organization of the paper is as follows. In section $2$ we will
give a necessary background and notation. In section $4$, using the
sequential convergence of the $\tau$-flow (that has been proved in
\cite{natasa2004}), we will construct a gauge on time intervals of an
arbitrary length, so that in the chosen gauge the $\tau$-flow equation
becomes strongly parabolic. In section $5$ we will use the
integrability assumption to prove that a soliton that we get in the
limit is unique, up to a diffeomorphism.

{\bf Acknowledgments:} I would like to thank my advisor Gang Tian for
bringing this problem to my attention, for many useful discussions,
suggestions and his constant support. Richard Hamilton, Huai-Dong Cao,
Tom Ilmanen, Peter Topping and Jeff Viaclovsky deserve many thanks as
well.
\end{section}

\begin{section}{Background}

Perelman's functional $\mathcal{W}$ and its properties will play an
important role in the paper. $M$ will always denote a closed
manifold. $\mathcal{W}$ has been introduced in \cite{perelman2002}.

$$\mathcal{W}(g,f,\tau) = (4\pi\tau)^{-\frac{n}{2}}\int_M
e^{-f}[\tau(|\nabla f|^2 + R) + f - n] dV_g.$$
We will consider this functional restricted to $f$ satisfying
\begin{equation}
\label{equation-equation_constraint}
\int_M (4\pi\tau)^{-\frac{n}{2}}e^{-f} dV = 1.
\end{equation} 
$\mathcal{W}$ is invariant under simultaneous scalings of $\tau$ and
$g$ and under a diffeomorphism change, i.e. $\mathcal{W}(g,f,\tau) =
\mathcal{W}(c\phi^*g,\phi^*f,c\tau)$ for a constant $c>0$ and a
diffeomorphism $\phi$. Perelman showed that the Ricci flow can be
viewed as a gradient flow of a functional $\mathcal{W}$, which is one
of the reasons why this functional plays an important role throughout
\cite{perelman2002}. Let $\mu(g,\tau) = \inf\mathcal{W}(g,f,\tau)$
over smooth $f$ satisfying (\ref{equation-equation_constraint}). It
has been showed by Perelman that $\mu(g,\tau)$ is achieved by some
smooth function $f$ on a closed manifold $M$, that $\mu(g,\tau)$ is
negative for small $\tau > 0$ and that it tends to zero as $\tau\to
0$.  

We will explain the motivation why we have decided to study this flow
instead of a normalized one in which a volume of a manifold has been
fixed along the flow. First of all, there is a simple
reparametrization that allows us to go from a $\tau$-flow to an
unnormalized flow and many smoothing regularity properties that have
been proved for the unnormalized flow continue to hold for a
$\tau$-flow as well. For example, Hamilton's compactness theorem also
holds for the $\tau$-flow. This is because Shi's estimates hold for
$\tau$-flow as well, and therefore, since we have a uniform curvature
bound on the solutions to a $\tau$-flow, we may assume uniform bounds
on all covariant derivatives of the curvature, $|D^p\rem| \le
C(p)$. The reparametrization that we use to go from a $\tau$-flow to
an unnormalized flow is as follows. Let $c(s)=1-\frac{s}{\tau}$ and
$t(s) = -\tau\ln(1-\frac{s}{\tau})$. Let $\tilde{g}(s) =
c(s)g(t(s))$. $\tilde{g}(s)$ is a solution to an unnormalized Ricci
flow. On the other hand we have that
$\mathcal{W}(g(t(s)),f(t(s)),\tau) =
\mathcal{W}(\tilde{g}(s),\tilde{f}(s),\tau - s)$. By the monotonicity
formula for $\mathcal{W}$ we have that the later quantity is
increasing along an unnormalized Ricci flow and therefore the former
quantity is increasing along the $\tau$ flow as well. The monotonicity
formula for a $\tau$-flow gets the simpler form;
$\mathcal{W}(g(t),f(t),\tau)$ is increasing along the $\tau$-flow,
while $f(t)$ changes by the evolution equation $\frac{d}{dt}f =
-\Delta f + |\nabla f|^2 -R + \frac{n}{2\tau}$ and $\tau$ is just a
constant. The fact that $\tau$ is now a constant will be very useful
in taking the limits of the minimizers for $\mathcal{W}$.

One of the most important properties of $\mathcal{W}$ is the
monotonicity formula.

\begin{theorem}[Perelman]
\label{theorem-theorem_monotonicity}
$\frac{d}{dt}\mathcal{W} = \int_M 2\tau|R_{ij} + \nabla_i\nabla_j f -
\frac{1}{2\tau}g_{ij}|^2(4\pi\tau)^{-\frac{n}{2}}e^{-f} dV \ge 0$ and
therefore $\mathcal{W}$ is increasing along the flow described by the
following equations

\begin{eqnarray*}
\frac{d}{dt}g_{ij} &=& -2R_{ij}, \\
\frac{d}{dt}f &=& -\Delta f + |\nabla f|^2 - R + \frac{n}{2\tau}, \\
\dot{\tau} &=& -1.
\end{eqnarray*}
\end{theorem}

One of the very important applications of the monotonicity formula is
noncollapsing theorem for the Ricci flow that has been proved by
Perelman in \cite{perelman2002}.

\begin{definition}
Let $g_{ij}(t)$ be a smooth solution to the Ricci flow
$(g_{ij})_t=-2R_{ij}(t)$ on $[0,T)$. We say that $g_{ij}(t)$ is
locally collapsing at $T$, if there is a sequence of times $t_k\to T$
and a sequence of metric balls $B_k = B(p_k,r_k)$ at times $t_k$, such
that $\frac{r_k2}{t_k}$ is bounded, $|\rem|(g_{ij}(t_k))\le r_k^{-2}$ in
$B_k$ and $r_k^{-n}\vol(B_k)\to 0$.
\end{definition}

\begin{theorem}[Perelman]
\label{theorem-perelman_theorem}
If $M$ is closed and $T < \infty$, then $g_{ij}(t)$ is not locally
collapsing at $T$.
\end{theorem}

The corollary of Theorem \ref{theorem-perelman_theorem} is 

\begin{corollary}
\label{corollary-corollary_ancient}
Let $g_{ij}(t)$, $t\in [0,T)$ be a solution to the Ricci flow on a
closed manifold $M$, where $T < \infty$. Assume that for some
sequences $t_k\to T$, $p_k\in M$ and some constant $C$ we have $Q_k =
|\rem|(x,t) \le CQ_k$, whenever $t < t_k$. Then a subsequence of
scalings of $g_{ij}(t_k)$ at $p_k$ with factors $Q_k$ converges to a
complete ancient solution to the Ricci flow, which is
$\kappa$-noncollapsed on all scales for some $\kappa > 0$.
\end{corollary}

We would like to recall a definition of a soliton that will appear in
later sections.

\begin{definition}
A Ricci soliton $g(t)$ is a solution to a Ricci flow equation that
moves by $1$-parameter group of diffeomorphisms $\phi(t)$, i.e. $g(t)
= \phi(t)^*g(0)$.
\end{definition}
The equation for a metric to move by a diffeomorphism in the direction
of a vector field $V$ is $2\ric(g) = \mathcal{L}_V(g)$, or $R_{ij} =
g_{ik}D_jV^k + g_{jk}D_iV^k$. If the vector field $V$ is the gradient
of a function $f$, we say that the soliton is the gradient Ricci
soliton. Moreover, we can consider the solutions to the Ricci flow
that move by diffeomorphisms and also shrink or expand by a factor at
the same time. The stationary solutions of the unnormalized Ricci flow
are the Ricci flat metrics. The Ricci solitons are the generalizations
of those, namely they are the stationary solutions to the Ricci flow
equations, up to diffeomorphisms. 

\end{section}

\begin{section}{Uniqueness of a limit soliton}

In \cite{natasa2004} we have proved the sequential convergence of a
$\tau$-flow with uniformly bounded curvatures and diameters toward the
solitons. In this section we will assume that one of the limit solitons
is integrable, in order to prove the uniqueness of a soliton in the
limit, up to a diffeomorphism. We will first construct a gauge in
which a $\tau$-flow becomes a strictly parabolic flow. Similar ideas
to those in \cite{cheeger1994} will help us finish the proof of
Theorem \ref{theorem-theorem_uniqueness}.

\begin{subsection}{The construction of a gauge}

To construct the right gauge, assume for simplicity that we are in a
situation when $g(t) \to h$ as $t\to\infty$, where $h$ is an Einstein
metric, with the Einstein constant $\frac{1}{2\tau}$. We will see how
we construct a gauge so that our modified Ricci flow equation becomes
strictly parabolic on time intervals of an arbitrary length, if we go
sufficiently far in time direction. This construction applies to our
more general case, just with minor modifications and only for
simplicity reasons we have decided to consider a case of an Einstein
metric in a limit. The main purpose of this section is to prove the
following Proposition that will be reformulated in the next section
for our more general setting.

\begin{proposition}
\label{proposition-proposition_extension}
Let $A>0$ be an arbitrary real number, $k$ an integer and $0 < \alpha
< 1$. There exists $\epsilon_0(A,k)$ such that for every $\epsilon <
\epsilon_0$ there exists $s_0 = s_0(A,k,h,\epsilon)$, such that for all
$t_0 \ge s_0$ the equation
\begin{eqnarray}
\label{equation-equation_gauge1}
\frac{d}{dt}\phi &=& \Delta_{g(t),h}\phi,\\
\phi(t_0) &=& \phi_{t_0}, \nonumber
\end{eqnarray}
has a solution $\phi(t)$, so that it is a diffeomorphism, $|\phi(t) -
\Id|_{k,\alpha,h} < \epsilon$ and $|\phi^*g(t) - h_0|_{k,\alpha} <
\epsilon$, for every $t\in[t_0,t_0+A]$. $\phi_{t_0}$ is chosen to be a
diffeomorphism so that $\delta_{\phi^*(t_0)h}(g(t_0)) = 0$.
\end{proposition}

\begin{definition}
Let $\phi: M\to M$ be a smooth function. Define $e(\phi-\Id) =
g^{ij}h_{kl}(\phi_i^k - \Id_i^k)(\phi_j^l - \Id_j^l)$. Define
$E(\phi-\Id) = \int_M e(\phi-\Id)$ and $F_l = \phi^l - \Id^l$.
\end{definition}

Throughout the proof of Proposition
\ref{proposition-proposition_extension} we will have a tendency to use
the same symbol for different uniform constants.

\begin{proof}[Proof of Proposition \ref{proposition-proposition_extension}]
Fix $A > 0$. Let $\epsilon > 0$ be very small (we will see later how
small we want to take it). We know that for $s_0$ sufficiently big we
can make $|g(t) - h|$ as small as we want, and therefore we have that
$\delta_{\phi(t_0)^*h}g(t_0) = 0$ implies that $|\phi(t_0) -
\Id|_{k+2,\alpha,h} < \epsilon/1000$ on $M$ (see \cite{cheeger1994}
for more details). Choose some $t_0 \ge s_0$. We can make
$|F(t_0)|_{N,\alpha,h}$, for say $N >> k$ as small as we want by
choosing $s_0$ sufficiently big. Since $g(t)\to h$ as $t\to\infty$,
the coefficients and the initial data of harmonic map flow
(\ref{equation-equation_gauge1}) are uniformly bounded and uniformly
close to each other for $t_0$ big enough. This implies that there
exists a uniform constant $\delta_1 > 0$ so that a solution to
(\ref{equation-equation_gauge1}) exists on $t\in[t_0,t_0+\delta_1)$,
for all $t_0 \ge s_0$. For the same reasons there exists some $\delta
> 0$ such that $|F(t)|_{W^{2,N},g(t)} < \epsilon$, for $t\in
[t_0,t_0+\delta)$. We can assume that we have chosen $N$ big enough so
that as a consequence of Sobolev embedding theorems we have that
$|F|_{k,\alpha, g(t)} < \tilde{\epsilon}$ ($\tilde{\epsilon}$ differs
from $\epsilon$ by a Sobolev embedding constant) for all $t\in
[t_0,t_0+\delta)$ and all $t_0 \ge s_0$. We want to show that the
estimate $|F(t)|_{W^{2,N},g(t)} < \epsilon $ holds past time
$t_0+\delta$, until $\delta < A$. Then $|F|_{k,\alpha,g(t)} <
\tilde{\epsilon}$ continues to hold past time $t_0+\delta$, until
$\delta < A$. This actually gives a uniform upper bound on the energy
densities on whole manifold $M$. To see this, notice that a bound
$|F|_{k,\alpha} < \tilde{\epsilon}$ implies that $e(\phi-\Id) \le
C\tilde{\epsilon}$. Since
$$e(\phi-\Id) = e(\phi) + e(\Id) - 2g^{ij}h_{kl}\Id_i^k\phi_j^l,$$ 
by the Schwartz inequality for quadratic forms and the interpolation
inequality we get that
\begin{eqnarray*}
e(\phi) &\le& C\tilde{\epsilon} + C + 2(g^{ij}h^{kl}\phi_i^k\phi_j^l)^{1/2}
(g^{ij}h^{kl}\Id_i^k\Id_j^l)^{1/2} \\
&\le& C\tilde{\epsilon} + C + \eta e(\phi),
\end{eqnarray*}
for some $\eta < 1$, which implies that $e(\phi) \le \tilde{C}$. By
the results proved by Eells and Sampson in \cite{eells1963} there
exists $\bar{\delta}$, depending on $(M,h)$ and the uniform bound on
the energy densities $\tilde{C}$, so that for every $s\in
[t_0,t_0+\delta)$ a solution to a harmonic map flow
(\ref{equation-equation_gauge1}) can be extended to
$[s,s+\bar{\delta}]$. If $t_0+\delta +\bar{\delta} < t_0 + A$, we can
repeat the procedure above for a solution $\phi(t)$, on time interval
$[t_0,t_0+\delta+\bar{\delta})$ to get that the energy density
estimates with the same constant $\tilde{C}$ hold past time
$t_0+\delta+\bar{\delta}$. Since all our estimates depend only on $A$
and the uniform bounds on geometries $g(t)$, we can iterate the
argument till we reach time $t_0+A$, for every $t_0\ge s_0$. As a
result, we will get $\phi(t)$, a solution to
(\ref{equation-equation_gauge1}), such that $|\phi(t) -
\Id|_{k,\alpha} < \tilde{\epsilon}$ for all $t\in [t_0,t_0+A]$.

We know that $(\Delta_{g(t),h} \Id)^{\gamma} =
g^{\alpha\beta}(\Gamma(h)_{\alpha\beta}^{\gamma} -
\Gamma(g)_{\alpha\beta}^{\gamma})$ and that $\frac{d}{dt}\Id = 0$.
Therefore, we have

\begin{equation}
\label{equation-equation_coordinate1}
\frac{d}{dt}(\phi^k - \Id^k) = \Delta_{g(t),h}(\phi^k - \Id^k) +
g^{ij}(\Gamma_{ij}^k(h) - \Gamma_{ij}^k(g)),
\end{equation} 
where we can choose $s_0$ so big, that the last term is 
arbitrarily small (since $g(t) \to h$). We will see later how 
small we want to make it, for now  we can say it is less 
than some $\epsilon_1 > 0$.

Before we start establishing the estimates on $F = \phi - \Id$, we
will occupy ourselves with the problem of replacing equation
(\ref{equation-equation_coordinate1}) which in terms of local
coordinates on $M$ is a local system of equations, by some much more
global system. Passing to a global system of equations will make
establishing the estimates on $F$ much easier. We will follow a
discussion in \cite{eells1963}.

Since $M$ is compact, there exists an embedding $\omega:M\to R^q$ and
due to Eells and Sampson (\cite{eells1963}) it is always possible to
construct a smooth Riemannian metric $g'' = (g''_{ab})_{1\le a,b\le
q}$ on a tubular neighborhood $N$ of $M$ in $R^q$, such that $N$ is
Riemannian fibered. They actually meant that if $\pi:N\to M$ is a
projection map, it suffices to construct an appropriate smooth inner
product in each space $R^q(p)$ for all $p\in M$, for which we can
translate that tangent space to any point $m\in N$ along the straight
line segment (that is contained in $N$) from $p = \pi(m)$ to
$m$. Following the arguments of section $7$ in \cite{eells1963} we
find that the evolution equation
(\ref{equation-equation_coordinate1}), given in local coordinates is
satisfied by $\phi - \Id$ if and only if $W - \tilde{W}$, where $W =
\omega\circ\phi$ and $\tilde{W} = \omega\circ\Id$ satisfies

\begin{equation}
\label{equation-equation_global1}
\frac{d}{dt}(W^c - \tilde{W}^c) = \Delta(W^c - \tilde{W}^c) +
\pi_{ab}^c(W_i^a - \tilde{W}_i^a)(W_j^b - \tilde{W}_j^b)g^{ij} +
\frac{\partial \omega^c}{\partial
y^k}g^{ij}(\Gamma_{ij}^k(h)-\Gamma_{ij}^k(g)),
\end{equation}
where $(y_1,\dots,y_n)$ are the local coordinates on $M$. Moreover,
since $M$ is compact, the projection $\pi$ satisfies (see
\cite{eells1963})
$$|\pi_{ab}^c|_{k+1,\alpha} \le C,$$
on $M$ and there are constants $A_1$ and $A_2$ so that
$$A_1ds_0^2 \le ds^2 \le A_2ds_0^2,$$ where $ds_0^2$ denotes the line
element induced on $M$ by the usual metric on $R^q$.  These estimates
immediately imply that
$$|\frac{\partial^k \pi_{ab}^c}{\partial y^k}W_i^aW_j^bg^{ij}| \le
C(k)e(\phi),$$ 
where also $e(\phi) = g''_{ab}W_i^aW_j^bg^{ij}$, 
$e(\phi - \Id) = g''_{ab}(W - \tilde{W})g^{ij}$.
Moreover, if $\tilde{F}^c = W^c - \tilde{W}^c$ then
$|\frac{\partial^k\pi_{ab}^c}{\partial
y^k}\tilde{F}_i^a\tilde{F}_j^bg^{ij}| \le Ce(\phi - \Id)$. 

The evolution equation for $e(\phi - \Id)$ (see for details
\cite{eells1963} and \cite{hamiltonMA}) is
\begin{eqnarray}
\label{equation-equation_schwartz}
\frac{d}{dt}e(\phi-\Id) &=& \Delta e(\phi-\Id) - 2|D^2(\phi -
\Id)|^2 + 2\rem(D(\phi-\Id),D(\phi-\Id),D(\phi-\Id),D(\phi-\Id)) \nonumber\\
&-& \frac{1}{\tau}e(\phi-\Id) + g^{ij}h_{kl}(\phi_j^k - \Id_j^k)[g^{pq}
(\Gamma_{pq}^l(h) - \Gamma_{pq}^l(g))]_i
\end{eqnarray}
where $\rem(D(\phi-\Id),D(\phi-\Id),D(\phi-\Id),D(\phi-\Id)) =
g^{ik}g^{jl}R_{pqmn}D_i(\phi^p - \Id^p)D_j(\phi^q - \Id_q) D_k(\phi^m
- \Id^m)D_l(\phi^n - \Id^n)$ and $|D^2(\phi - \Id)|^2 =
g^{ik}g^{jl}h_{pq}D^2_{ij}(\phi^p - \Id^p)D^2_{kl}(\phi^q - \Id^q)$.
Applying the Schwarz inequality for quadratic forms and using the fact
that $2\sqrt{\tau}(g^{pq}(\Gamma_{pq}^l(h) - \Gamma_{pq}^l(g)))_i$ can be made
arbitrarily small by choosing $s_0$ sufficiently big (e.g. smaller
than $\frac{2\epsilon}{1000}$), the last term in inequality
(\ref{equation-equation_schwartz}) can be estimated as

$$g^{ij}h_{kl}(\phi_j^k - \Id_j^k)[g^{pq} (\Gamma_{pq}^l(h) -
\Gamma_{pq}^l(g))]_i \le
\frac{e(\phi-\Id)^{\frac{1}{2}}}{2\sqrt{\tau}}(2\epsilon)/1000.$$ 
Factor of $1000$ (that we can increase if necessary) is chosen so that after
multiplying $\frac{\epsilon}{1000}$ by at most a polynomial expression
in $A$ (which will become more apparent later in the proof of
Proposition \ref{proposition-proposition_extension}) can be made again
much smaller than $\epsilon$. Therefore, for $t\in [t_0,t_0+\delta)$
we have that

\begin{claim}
\label{claim-claim_en_est1}
There exists $C$, small $\epsilon$ and sufficiently big $s_0$ such
that for all $t_0\ge s_0$ 
\begin{enumerate}
\item
$e(\phi - \Id) < \epsilon_1$,
\item
$E(\phi-\Id)(s) < \epsilon_1,$
\end{enumerate} 
for all $s$ belonging to a time
interval starting at $t_0$ at which $\phi$ exists, where $\epsilon_1$
is a constant that can be made much smaller than $\epsilon$.
\end{claim}

\begin{proof}
By using the interpolation inequality in
(\ref{equation-equation_schwartz}), we get 
\begin{eqnarray*}
\frac{d}{dt}e(\phi-\Id) &\le& \Delta e(\phi-\Id) + C\epsilon^4 -
\frac{1}{\tau}e(\phi-\Id) + \frac{1}{2\tau}e(\phi-\Id) + 
C\frac{\epsilon^2}{1000^2} \\
&\le& \Delta e(\phi-\Id) -\frac{1}{2\tau}e(\phi-\Id) + \frac{\epsilon}{1000},
\end{eqnarray*}
since we can start with $\epsilon$ as small as we want, in particular
we may choose $\epsilon$ so that $C\epsilon^4 + C\frac{\epsilon^2}{1000^2}
< \frac{\epsilon}{1000}$ and increase $s_0$ if necessary.

Let $f(t) = \max_Me(\phi-\Id)(t)$. Then
$$\frac{d}{dt}f \le -\frac{1}{2\tau}f + \frac{\epsilon}{1000},$$
$$\frac{d}{dt}f \le -\frac{1}{2\tau}(f - \frac{\tau\epsilon}{500}).$$
If we choose $s_0$ big enough, we may assume that $f(t_0) <
\frac{\tau\epsilon}{500}\vol_h(M)$. If $f(t) \ge
\frac{\tau\epsilon}{500}\vol_h(M))$ for some $t > t_0$, then $f(t)$ is
nonincreasing (because $\frac{d}{dt}f(t) \le 0$ and since it starts
as $f \le \frac{\tau\epsilon}{500}\vol_h(M))$), it will remain so
forever while $\phi$ exists. Denote by $\epsilon_1 =
\frac{\tau\epsilon}{500}\max_t\vol_{g(t)}(M)$.

\begin{equation}
\label{equation-equation_energy_estimate}
E(\phi-\Id)(s) = \int_U e(\phi-\Id)(s)dV_{g(s)} < \epsilon_1.
\end{equation}
\end{proof}

By Claim \ref{claim-claim_en_est1}, $e(W-\tilde{W})$ can be made much
smaller than $\epsilon$ whenever $\phi$ is defined (if $t_0$ is big
enough and $\epsilon$ is small enough).  The conditions $|F|_{W^{2,N}}
< \epsilon$ and $|F|_{k,\alpha} < \tilde{\epsilon}$ actually mean that
for $\tilde{F}$ we make an assumption that $|\tilde{F}|_{W^{2,N}} <
\epsilon$ and $|\tilde{F}|_{k,\alpha} < \tilde{\epsilon}$, for $t\in
[t_0,t_0+\delta)$ (these $\epsilon$ and $\tilde{\epsilon}$ can be
slightly different from those for $F$). In order to finish the proof
of Proposition \ref{proposition-proposition_extension} it is enough to
show that $|\tilde{F}|_{k,\alpha} < \tilde{\epsilon}$ continues to
hold past time $t_0+\delta$, for $t_0$ big enough.  From now on we
will consider a globally defined evolution equation
\begin{equation}
\label{equation-equation_global}
\frac{d}{dt}(\tilde{F}^c) = \Delta\tilde{F}^c +
\pi_{ab}^c\tilde{F}^a_i\tilde{F}^b_jg^{ij} + \frac{\partial
\omega^c}{\partial y^k}g^{ij}(\Gamma_{ij}^k(h)-\Gamma_{ij}^k(g)),
\end{equation}

\begin{step}
\label{step-step_1}
$\int_M |\tilde{F}^c|^2 dV_{g(t)}$ and
$\int_{t_0}^{t_0+\delta}\int_M|\nabla\tilde{F}^c|^2dV_{g(t)}$ can be
made much smaller than $\epsilon$, for all $t\in [t_0,t_0+\delta)$ and
for all $t_0$ big enough.
\end{step}

Multiply the equation (\ref{equation-equation_global}) by
$\tilde{F}^c$ and integrate it over $M$ against the metric
$g(t)$.

\begin{eqnarray}
\label{equation_equation_help_help}
\frac{1}{2}\frac{d}{dt}\int(\tilde{F}^c)^2 dV_{g(t)} &<& \int_M
(\tilde{F}^c)^2(\frac{n}{2\tau}-R)dV_{g(t)} - \int_M
|\nabla\tilde{F}^c|^2 d_{g(t)} +
\epsilon_1\int_M |\tilde{F}^c| dV_{g(t)} + \nonumber\\
&+& C(\int_M e(\tilde{F})^2
dV_{g(t)})^{1/2}(\int_M\tilde{F}^c)^2dV_{g(t)})^{1/2} \nonumber\\
&\le& \epsilon_1\epsilon  - \int_M|\nabla \tilde{F}^c|^2
\nonumber + \epsilon_1[\int_M |\tilde{F}^c| dV_{g(t)}] + 
C\epsilon\epsilon_1\int_M (\tilde{F}^c)^2dV_{g(t)}\nonumber \\ 
\end{eqnarray} 
since 
$$\int_M \tilde{F}^c\frac{\partial\omega^c}{\partial
y_l}g^{ij}(\Gamma(h_{ij}^l-\Gamma(g)_{ij}^l)dV_{g(t)} \le
\epsilon_1\int_M|\tilde{F}^c|dV_{g(t)},$$
\begin{eqnarray*}
\int_M\tilde{F}^cg^{ij}\pi_{ab}^c\tilde{F}_i^a\tilde{F}^b_j &\le& 
C(\int_M e(\tilde{F})^2)^{1/2}(\int_M(\tilde{F}^c)^2)^{1/2} \\
&<& C\epsilon\epsilon_1[\int(\tilde{F}^c)^2)]^{1/2}.
\end{eqnarray*} 
In the above estimates we have used the energy estimates
(\ref{equation-equation_energy_estimate}), the fact that $g(t)\to h$
as $t\to\infty$ uniformly on $M$ and that $|\tilde{F}|_{W^{2,N}} <
\epsilon$ for $t\in [t_0,t_0+\delta)$ (which implies
$|\tilde{F}|_{C^{k,\alpha}} < \tilde{\epsilon}$ for sufficiently big
$N$). For those reasons, $\epsilon_1 << \epsilon$ is a constant that
can be made much smaller than $\epsilon$, by taking $\epsilon$ small
and $s_0$ big. Integrate (\ref{equation_equation_help_help}) in $t$.
\begin{eqnarray*}
\frac{1}{2}\sup_{t\in[t_0,t_0+\delta)}\int(\tilde{F}^c)^2(t)dV_{g(t)} +
\sup_{t\in[t_0,t_0+\delta)}\int_{t_0}^{t}\int_M|\nabla\tilde{F}^c|^2dV_{g(t)} 
&\le& \int \frac{1}{2}(\tilde{F}^c)^2(t_0)dV_h + CA\epsilon\epsilon_1. 
\end{eqnarray*}
Since for big $t_0$ the first integral on the right hand side of the
previous inequality can be made much smaller than $\epsilon$, it
follows that for big $t_0$ and small $\epsilon$,

$$\sup_{t\in[t_0,t_0+\delta)}\int_M(\tilde{F}^c)^2(t)dV_{g(t)} < 
\tilde{\epsilon},$$
$$\sup_{t\in[t_0,t_0+\delta)}\int_M |\nabla\tilde{F}^c|^2dV_{g(t)} <
\tilde{\epsilon},$$
for some constant $\tilde{\epsilon} << \epsilon$ and these estimates
depend on $A$. 

\begin{step}
\label{step-step_3}
$\sup_{t\in[t_0,t_0+\delta)}\int_{t_0}^t\int_M|\frac{d}{dt}\tilde {F}^c|^2
dV_{g(t)}$ and $\sup_{t\in[t_0,t_0+\delta)}\int_{t_0}^t\int_M|\nabla^2
\tilde{F}^c|^2dV_{g(t)}$ can be made much smaller than $\epsilon$ for
big enough $s_0$ which depends on $A$ and on the rate of convergence of
$g(t)$ to $h$, for small enough $\epsilon$.
\end{step}

$\frac{d}{dt}\tilde{F}^c = \Delta \tilde{F}^c + H^c$, where $H^c =
\frac{\partial \omega^c}{\partial y_l}g^{ij}
(\Gamma(h)_{ij}^l-\Gamma(g)_{ij}^l) + g^{ij}\pi_{ab}^c\tilde{F}^a_i
\tilde{F}^b_j$. Then

\begin{equation}
\label{equation-equation_main}
(H^c)^2 = (\Delta \tilde{F}^c)^2 + (\frac{d}{dt}\tilde{F}^c)^2 - 2\Delta
\tilde{F}^c\frac{d}{dt}\tilde{F}^c.
\end{equation}

\begin{eqnarray}
\label{equation-equation_comb1}
-\int_M\frac{d}{dt}\tilde{F}^c\Delta \tilde{F}^c &=& \int_M
g^{ij}\nabla_i(\frac{d}{dt}\tilde{F}^c)\nabla_j \tilde{F}^c \\
&=& \frac{1}{2}\frac{d}{dt}\int_M |\nabla \tilde{F}^c|^2 +
\frac{1}{2}\int_M g^{ip}g^{jq}(-2R_{pq}+\frac{1}{\tau}g_{pq})
|\nabla \tilde{F}^c|^2 \nonumber \\ 
&-& \frac{1}{2}\int_M |\nabla \tilde{F}^c|^2(\frac{n}{2\tau}-R). \nonumber
\end{eqnarray}

\begin{eqnarray}
\label{equation-equation_comb2}
\int_M(\Delta \tilde{F}^c)^2 &=& \int_M|\nabla^2 \tilde{F}^c|^2 + 
\int_M g^{ij}g^{ks}\nabla_j \tilde{F}^cR_{jkp}^s\nabla_p \tilde{F}^c.
\end{eqnarray}
Combining (\ref{equation-equation_main}), (\ref{equation-equation_comb1})
and (\ref{equation-equation_comb2}) we get

\begin{eqnarray}
\label{equation-equation_comb}
& &\int_M|\frac{d}{dt}\tilde{F}^c|^2 + \int_M|\nabla^2\tilde{F}^c|^2 +
\frac{d}{dt}\int_M|\nabla \tilde{F}^c|^2 \le \nonumber \\
&\le& \int_M(H^c)^2 + C\int_M|\nabla\tilde{F}^c|^2dV_{g(t)}.
\end{eqnarray}
since
$$\int_M g^{ip}g^{jq}(2R_{pq}-\frac{1}{\tau}g_{pq})|\nabla \tilde{F}^c|^2
\le \epsilon_1\int|\nabla \tilde{F}^c|^2,$$
$$\int_M g^{ij}g^{ks}\nabla_j \tilde{F}^cR_{jkp}^s\nabla_p \tilde{F}^c \le 
C\int_M|\nabla \tilde{F}^c|^2,$$
Notice also that
\begin{eqnarray}
\label{equation-equation_H_c}
\int_M(H^c)^2dV_{g(t)} &\le& \int_M
C(g^{ij}(\Gamma_{ij}^l(h)-\Gamma_{ij}^l(g)))^2 dV_{g(t)} + C\int_M
e(\tilde{F})^2dV_{g(t)} \nonumber \\
&<& \epsilon_1 + C\epsilon\epsilon_1,
\end{eqnarray} 
since $e(\tilde{F}) < \epsilon$ for $t\in [t_0,t_0+\delta)$ and
$\int_M e(\tilde{F})dV_{g(t)} < \epsilon_1$ by energy estimates
(\ref{equation-equation_energy_estimate}). We will sometimes use the
same constant $\epsilon_1$ to denote any constant that can be made
much smaller than $\epsilon$, (the estimates above are possible if we
start with $\epsilon$ small enough and increase $s_0$ if necessary,
depending on how big $A$ is). If we integrate
(\ref{equation-equation_comb}) in $t$ and use the above estimates, we
get

\begin{eqnarray}
\label{equation-equation_together}
& &\int_M|\nabla \tilde{F}^c(t)|^2 dV_{g(t)} + 
\int_{t_0}^t\int_M|\frac{d}{dt}\tilde{F}^c|^2
+ \int_{t_0}^t\int_M|\nabla^2\tilde{F}^c|^2 \le \\
&\le& \int_M|\nabla \tilde{F}^c(t_0)|^2 dV_{g(t_0)} + 
\int_{t_0}^t\int_M(H^c)^2 + C\int_{t_0}^{t}\int_M|\nabla\tilde{F}^c|^2dV_{g(s)}ds
\nonumber \\
&\le& \epsilon_1 \nonumber,
\end{eqnarray} 
because of Step \ref{step-step_1}, the fact that $\int_M|H^c|^2 <
\epsilon_1$ for big $t_0$ and the fact that for big $t_0$ 
$\int_M|\nabla^k \tilde{F}|^2(t_0)dV_{g(t_0)}$ can be made very small.

\begin{step}
\label{step-step_4}
$\sup_{t\in[t_0,t_0+\delta)}\int_M|\frac{d}{ds}\tilde{F}^c|^2$,
$\int_{t_0}^t\int_M |\nabla\frac{d}{dt}\tilde{F}^c|^2$ can be made much smaller
than $\epsilon$ for big $t_0$, for all $t\in [t_0,t_0+\delta)$.
\end{step}
 
Let $\hat{F}^c = \frac{d}{dt}\tilde{F}^c$. Then
$$\frac{d}{dt}\hat{F}^c = \Delta\hat{F}^c +
g^{ip}g^{jq}(2R_{pq}-\frac{1}{\tau}g_{pq})\nabla_i\nabla_j \tilde{F}^c + 
\frac{d}{dt}H^c - g^{ij}\frac{d}{dt}\Gamma_{ij}^k\nabla_k\tilde{F}^c.$$
Multiply this equation by $\hat{F}^c$ and integrate it over $M$.
\begin{eqnarray*}
\frac{1}{2}\frac{d}{dt}\int_M(\hat{F}^c)^2 - \frac{1}{2}\int_M
(\hat{F}^c)^2(\frac{n}{2\tau}-R) &=&
-\int_M |\nabla \hat{F}^c|^2 +
\int_M g^{ip}g^{jq}(2R_{pq}-
\frac{1}{\tau}g_{pq})\nabla_i\nabla_j \tilde{F}^c\hat{F}^c  \\
&+& \int_M\hat{F}^c\frac{d}{dt}H^c - 
\int_M g^{ij}\frac{d}{dt}\Gamma_{ij}^k\nabla_k\tilde{F}^c\hat{F}^c.
\end{eqnarray*}
Integrate it in $t$ to get

\begin{eqnarray}
\label{equation-equation_again}
\frac{1}{2}\int_M(\hat{F}^c(t))^2 + \int_{t_0}^t\int_M|\nabla\hat{F}^c|^2
&\le& \frac{1}{2}\int_M(\hat{F}^c(t_0))^2 + \frac{1}{2}\int_{t_0}^t
\int_M(\hat{F}^c)^2(\frac{n}{2\tau}-R) + \nonumber \\
&+&\epsilon_1(\int_{t_0}^t\int |\nabla^2\tilde{F}^c|^2)^{1/2}(\int_{t_0}^t\int
|\hat{F}^c|^2)^{1/2}  \nonumber \\
&+& (\int_{t_0}^t\int_M (\hat{F}^c)^2)^{1/2}
(\int_{t_0}^t\int_M\frac{d}{dt}H^c)^2)^{1/2} \nonumber \\
&+& C(\int_{t_0}^t\int_M|\nabla \tilde{F}^c|^2)^{1/2}(\int_{t_0}^t
\int_M(\hat{F}^c)^2)^{1/2}.
\end{eqnarray}
Notice that
\begin{eqnarray}
\label{equation-equation_H}
\int_{t_0}^t\int_M(\frac{d}{dt}H^c)^2 &\le& C(\int_{t_0}^t\int_M
(\frac{d}{dt}(\frac{\partial\omega^c}{\partial y^l}g^{ij}
(\Gamma(h)_{ij}^l-\Gamma(g)_{ij}^l)))^2 \nonumber \\
&+& \int_{t_0}^t\int_M(((\frac{d}{dt}g^{ij})\pi_{ab}^c
\tilde{F}_i^a\tilde{F}_j^b)^2 \\
&+& \int_{t_0}^t\int_M(g^{ij}\pi_{ab}^c(\frac{d}{dt}\tilde{F}_i^a)
\tilde{F}_j^b)^2),
\end{eqnarray}
where
$$\int_{t_0}^t\int_M (\frac{d}{dt}(g^{ij}\frac{\partial\omega^c}{\partial
y^l}(\Gamma(h)_{ij}^l-\Gamma(g)_{ij}^l)))^2 < \epsilon_1,$$
$$\int_{t_0}^t\int_M((\frac{d}{dt}g^{ij})\pi_{ab}^c\tilde{F}_i^a
\tilde{F}_j^b)^2 \le \epsilon_1,$$
if $t_0$ is big enough, since $g(t)\to h$ uniformly on $M$ and
$\frac{d}{dt}g^{ij} = g^{pi}g^{qj}(2R_{pq}-\frac{1}{\tau}g_{pq})$, and
$$\int_{t_0}^t\int_M(g^{ij}\pi_{ab}^c
(\frac{d}{dt}\tilde{F}_i^a)\tilde{F}_j^b)^2
\le C\epsilon\int_{t_0}^t\int_M|\nabla\hat{F}^c|^2 < 
\frac{1}{2}\int_{t_0}^t\int_M|\nabla\hat{F}^c|^2,$$
if we choose $\epsilon$ small enough, such that $C\epsilon < \frac{1}{2}$,
since $|\nabla_j \tilde{F}| < \epsilon$ for $t\in [t_0,t_0+\delta)$.
$$\int_M|\hat{F}^c|^2dV_{g(t)} \le C(\int_{t_0}^t\int_M|\nabla^2\tilde{F}^c|^2
+ C\int_{t_0}^{t_0+\delta}\int_M|\nabla\tilde{F}^c|^2 + \epsilon_1) \le
C(\epsilon_1 + \epsilon_1),$$
by Step \ref{step-step_3}.
The assertion of Step \ref{step-step_4} follows now immediately from
(\ref{equation-equation_again}).

From the estimate (\ref{equation-equation_H}) we can now get (using the
estimates of Steps \ref{step-step_1}, \ref{step-step_3} and
\ref{step-step_4}) that $\int_{t_0}^t\int_U(\frac{d}{dt}H^c)^2$ can be
much smaller than $\epsilon$. Consider the equation
\begin{equation}
\label{equation-equation_hat}
\frac{d}{dt}\hat{F}^c = \Delta\hat{F}^c + \hat{H}^c,
\end{equation}
where $\hat{F}^c = \frac{d}{dt}\tilde{F}^c$ and $\hat{H}^c =
g^{ip}g^{jq}(2R_{pq}-\frac{1}{\tau}g_{pq})\nabla_i\nabla_j \tilde{F}^c +
\frac{d}{dt}H^c -
g^{ij}\nabla_k\tilde{F}^c\frac{d}{dt}(\Gamma_{ij}^k)$. Since $g(t)\to
h$, where $\ric(h) = \frac{1}{2\tau}h$, by using the previous
estimates, we can easily see that $\int_{t_0}^t\int_U (\hat{H}^c)^2
dV_{g(s)}ds$ can be made much smaller than $\epsilon$. In the same
manner as we have obtained the estimates in step \ref{step-step_3} for
$\tilde{F}^c$, we can get the following estimates for $\hat{F}^c =
\frac{d}{dt}\tilde{F}^c$ by considering the evolution equation
(\ref{equation-equation_hat}).
$$\sup_{t\in[t_0,t_0+\delta)}\int_U|\nabla\frac{d}{dt}\tilde{F}^c|^2,$$
$$\int_{t_0}^{t_0+\delta}\int_U|\nabla^2\frac{d}{dt}\tilde{F}^c|^2,$$
$$\int_{t_0}^{t_0+\delta}\int_U(\frac{d^2}{dt^2}\tilde{F}^c)^2,$$
can be made much smaller than $\epsilon$ for big $t_0$.

We have that $\Delta \tilde{F}^c = \frac{d}{dt}\tilde{F}^c - H^c$
where $W^{1,2}$ norm of $\frac{d}{dt}\tilde{F}^c$ and $L^2$ norm of
$H^c$ can be made much smaller than $\epsilon$. By elliptic regularity
theory we can get that $W^{2,2}$ norm of $\tilde{F}^c$ can be made
much smaller than $\epsilon$ (since it can be estimated in terms of
$W^{1,2}$ norm of $\frac{d}{dt}\tilde{F}^c$ and $L^2$ norm of
$H^c$). Using that and the fact that $|\tilde{F}|_{1,\alpha,g} <
\epsilon$ notice that
\begin{eqnarray*}
\int_M|\nabla H^c|^2 &<& C(\epsilon_1 + \int_M g^{rs}
(\nabla_r(g^{ij}\pi_{ab}^c))\tilde{F}_i^a\tilde{F}_j^b)
(\nabla_sg^{i'j'}\pi_{a'b'}^{c'}\tilde{F}_{i'}^{a'}\tilde{F}_{j'}^{b'}) \\ 
&+& \int g^{rs}(g^{ij}\pi_{ab}^c\nabla_r \tilde{F}_i^a \tilde{F}_j^b)
(g^{i'j'}\pi_{a'b'}^{c'}\nabla_s\tilde{F}_{i'}^{a'}\tilde{F}_{j'}^{b'}) \\ 
&<& C(\epsilon_1 + \epsilon\sum_a\int_M|\nabla\tilde{F}^a|^2 + 
\epsilon\sum_b\int_M |\nabla^2 \tilde{F}^b|^2) \\
&<& \tilde{\epsilon},
\end{eqnarray*}
for some small constant $\tilde{\epsilon}$, that can be assumed to be
much smaller than $\epsilon$, since $W^{2,2}$ norm of $\tilde{F}^c$
can be made much smaller than $\epsilon$. By elliptic regularity
theory this implies that $W^{3,2}$ norm of $\tilde{F}^c$ can be made
much smaller than $\epsilon$ for $t_0$ very big.

We can continue our proof by studying the equation
$\frac{d}{dt}\hat{F}^c = \Delta \hat{F}^c +
\hat{H}^c$. $|\tilde{F}^c|_{2,\alpha} < \epsilon$ for $t\in
[t_0,t_0+\delta)$. By a standard parabolic regularity we can get the
higher order estimates of $\tilde{F}^c$, by constants that are
comparable to $\epsilon$. Therefore, by the similar analysis as above
we can get that $W^{3,2}$ norms of $\hat{F}^c$ can be made much
smaller than $\epsilon$, since from the estimates that we have got
till this point we can again easily get that $W^{1,2}$ norm of
$\hat{H}^c$ can be made much smaller than $\epsilon$. Consider again
the equation
\begin{equation}
\label{equation-equation_ponovo}
\Delta \tilde{F}^c = \frac{d}{dt}\tilde{F}^c - H^c.
\end{equation} 
We know that $W^{3,2}$ norm of $\frac{d}{dt}\tilde{F}^c$ and $W^{3,2}$ norms of
$\tilde{F}^c$ can be made much smaller than $\epsilon$. Let's check
that $W^{3,2}$ norm of $H^c$ can be made much smaller than $\epsilon$
as well. In order for it to be true it is enough to check that $\int_M
|\nabla^3 (g^{ij}\pi_{ab}^c\tilde{F}_i^a\tilde{F}_j^b)|^2$ can be made much
smaller than $\epsilon$.
\begin{eqnarray*}
\int |\nabla^3 (g^{ij}\pi_{ab}^c\tilde{F}_i^a\tilde{F}_j^b)|^2 &\le&
C( \epsilon\int_M|\nabla\tilde{F}|^2 + C\sum_{a,b}\int
|\nabla^4\tilde{F}^a|^2|\nabla \tilde{F}^b|^2 + C\sum_{a,b}\int_M
|\nabla^3\tilde{F}^a|^2|\nabla^2 \tilde{F}^b|^2) \\ 
&<& C(\epsilon_1 + C\epsilon\epsilon_1 + \epsilon_1^2) < \tilde{\epsilon},
\end{eqnarray*}
since $|\tilde{F}|_{W^{2,N},g(t)} < \epsilon$ for all $t\in[t_0,t_0+\delta)$ and all
$t_0\ge s_0$. From here, again by elliptic regularity theory applied 
to equation (\ref{equation-equation_ponovo}), it
follows that $W^{5,2}$ norm of $\tilde{F}^c$ can be made much smaller than
$\epsilon$.

We can continue the proof in a similar manner as above, by taking the
higher order derivatives of our original equation $\Delta \tilde{F}^c
= \frac{d}{dt}\tilde{F}^c - H^c$ in $t$, using the estimates that we
get on the way and then go backward to our original equation to
improve a regularity of $\tilde{F}^c$. As a result, we can get (performing the
previously described procedure sufficiently many times) that
$|\tilde{F}|_{W^{N,2},g(t)} < \epsilon$ continues to hold past time
$t_0+\delta$.
\end{proof}

So far we have proved that for every $A>0$ and an integer $k$ there
exists $\epsilon_0=\epsilon_0(A,k)$ such that for every $\epsilon <
\epsilon_0$ we can find $s_0 = s_0(A,\epsilon,k)$, so that $\forall \:\:
t_0 \ge s_0$ there exists a solution of
\begin{eqnarray}
\label{equation-equation_gauge}
\frac{d}{dt}\phi(t) &=& \Delta_{g(t),h}\phi(t)\\
\phi(t_0) &=& \phi_{t_0} \nonumber,
\end{eqnarray}
for all $t\in[t_0,t_0+A]$ and $|\phi - \Id|_{k,\alpha} < \epsilon$.
 
We want to show that these maps $\phi(t):M \to M$ are actually 
diffeomorphisms which will imply that we have constructed an $1$-
parameter family of gauges such that for 
$\bar{g}(t) = (\phi(t)^*)^{-1}g(t)$ the linearization of the Ricci-DeTurck
flow
$$\frac{d}{dt}\bar{g} = -2\bar{R}_{ij} + \frac{1}{\tau}\bar{g}_{ij} +
\nabla_i W_j + \nabla_j W_i,$$ 
with $\bar{g}(t_0) =
(\phi_{t_0}^{-1})^*g(t_0)$ is strictly parabolic ($W_j =
\bar{g}_{jk}\bar{g}^{pq}(\Gamma_{pq}^k(\bar{g}) - \Gamma(h)_{pq}^k)$).

\begin{corollary}
Adopt the notation from Proposition
\ref{proposition-proposition_extension}. $\phi(t)$ are
diffeomorphisms for all $t\in [t_0,t_0+A]$ and all $t_0\ge s_0$.
\end{corollary}

\begin{proof}
Fix any $t_0 \ge s_0$. Consider the equation
\begin{eqnarray}
\label{equation-equation_parabolic}
\frac{d}{dt}\tilde{g}_{ij} &=& -2\tilde{R}_{ij} +
\frac{1}{\tau}\tilde{g}_{ij} + \nabla_i V^j + \nabla_j V^i,\\
\tilde{g}(t_0) &=& (\phi_{t_0}^{-1})^*g(t_0) \nonumber,
\end{eqnarray}
where $V^k = \tilde{g}^{pq}(\Gamma^k_{pq}(\tilde{g}) - \Gamma^k_{pq}(h))$.
This is a strictly parabolic system of equations and therefore
there exists some $\delta > 0$ so that a solution $\tilde{g}$
exists for all times $t\in [t_0,t_0+\delta)$. On the other hand, look
at the system
\begin{eqnarray}
\label{equation-equation_ODE}
\frac{d}{dt}\psi(t) &=& -V\circ\psi(t),\\
\psi(t_0) &=& \phi_{t_0} \nonumber.
\end{eqnarray}
Vector fields $V(t)$ are defined for $t\in [t_0,t_0+\delta)$ and
therefore the system (\ref{equation-equation_ODE}) has a solution
$\psi(t)$ for all those times. It is easy to show (a classical result)
that all $\psi(t)$ are diffeomorphisms for $t\in [t_0,t_0+\delta)$. The
simple computation (due to the fact that $g(t)$ is a solution of the
Ricci flow equation) shows that
$$\frac{d}{dt}\psi(t) = \Delta_{g(t),h}\psi(t),$$ with $\psi(t_0) =
\phi_{t_0}$. Because of the uniqueness of a harmonic map flow with the
same initial data (we know that our solutions are smooth and uniformly
bounded, so the uniqueness follows by the arguments of Eells and
Sampson in \cite{eells1963}), we have that $\psi(t) = \phi(t)$ for all
$t\in [t_0,t_0+\delta)$. This means $\phi(t)$ is a diffeomorphism for
$t\in [t_0,t_0+\delta)$ and $\tilde{g}(t) = (\phi(t)^{-1})^*g(t)$. We
know that for all $t\in [t_0,t_0+A]$, for $t_0$ sufficiently big, we
have that $|\phi(t) - \Id|_{k,\alpha,h} < \epsilon$. Therefore,
$|\phi^{-1} - \Id|_{k,\alpha}$ can be made small which implies that
$|\tilde{g}(t) - g(t)|$ can be made very small, comparable to
$\epsilon$, for all $t\in [t_0,t_0+\delta)$.  We want to extend a
solution $\tilde{g}(t)$ of (\ref{equation-equation_parabolic}) all the
way up to $t_0+A$. Since $|\tilde{g}(t) - h| < \tilde{\epsilon}$ and
since our flow (\ref{equation-equation_parabolic}) is strictly
parabolic, there exists $t_1 = t_1(h,\tilde{\epsilon})$ so that for
every $t\in [t_0,t_0+\delta)$, a solution to
(\ref{equation-equation_parabolic}) exists for all times $s\in
[t,t+t_1)$. That means we can extend our solution past time
$t_0+\delta$. Since our estimates on $|\tilde{g}(t) - h|$ for those
times for which a solution $\tilde{g}(t)$ exists are independent of
$\delta \le A$, we can easily extend our solution all the way up to
$t_0+A$, with $|\tilde{g}(t) - h|$ staying very small (comparable to
$\epsilon$) for all $t\in [t_0,t_0+A]$. Existence of $\tilde{g}(t)$
for $t\in[t_0,t_0+A]$ gives that $\phi(t)$ stays a diffeomorphism for
all times up to $t_0+A$, because it solves the equation
(\ref{equation-equation_ODE}).
\end{proof}

\end{subsection}

\begin{subsection}{The integrable case}

The proofs in this subsection are motivated by those in
\cite{cheeger1994}, where Cheeger and Tian have considered the
uniqueness problem of tangent cones under the assumption of
integrability of one of the tangent cones and under some curvature and
volume bounds.

\begin{remark}
\label{remark-remark_proved}
So far we have proved that for every $A>0$ and an integer $k$ there
exists $\epsilon_0 = \epsilon_0(A,k)$ such that for every $\epsilon
\le \epsilon_0$ there exists $s_0 = s_0(\epsilon,A,k)$ with the
property that for every $t_0\ge s_0$ there is an $1$-parameter family
of diffeomorphisms $\phi(t)$ so that
\begin{enumerate}
\item
$\phi^{-1}$ solves a harmonic map flow equation
\begin{eqnarray*}
\frac{d}{dt}\phi^{-1} &=& \Delta_{g,h}\phi^{-1},\\
\phi^{-1}(t_0) &=& \phi_{t_0},
\end{eqnarray*}
where $\delta_{\phi_{t_0}^*h}(g(t_0)) = 0$, for $t\in [t_0,t_0+A]$,
\item
$\tilde{g}=\phi^*g$ solves strictly parabolic equation on $[t_0,t_0+A]$
$$\frac{d}{dt}\tilde{g} = -2\ric(\tilde{g}) + \frac{1}{\tau}\tilde{g}
+ \nabla_iV_j + \nabla_jV_i,$$ 
where $V^i = \tilde{g}^{pq}(\Gamma^{i}_{pq}(\tilde{g}) - \Gamma^i_{pq}(h))$. 
We will say that $\tilde{g}$ is in a {\bf standard form} around $h$. We
will denote by $P_{h_0}(\tilde{g}) = \nabla_iV_j + \nabla_jV_i$.
\item
$|\phi - \Id|_{k,\alpha} < \epsilon$.
\item
$|\tilde{g} - h|_{k,\alpha} < \epsilon$.
\end{enumerate}
\end{remark}

From now on, we will simply write $\phi g$ instead of $\phi^*g$. By
the assumptions of Theorem \ref{theorem-theorem_uniqueness} there
exists a limit soliton, say $h(t)$ which is integrable. There is a
sequence $t_i$ such that $g(t_i+t)\to h(t)$ as $i\to\infty$ and
$$R_{ij} + \nabla_i\nabla_j f - \frac{1}{2\tau}h_{ij}(t) = 0,$$ 
for some function $f$. From before we know that $f(t)$ is a minimizer for
$\mathcal{W}$ with respect to a metric $h(t)$, for every $t$. Let
$\psi(t)$ be $1$-parameter family of diffeomorphisms induced by a
vector field $-\nabla f$. Then $h(t) = \psi^*(t)h_0$, where $h_0 =
h(0)$. Since $h_0 = (\psi^{-1})^*h(t)$, it satisfies the equation
\begin{equation}
\label{equation-equation_soliton_h0}
0 = \frac{d}{dt}h_0 = -2\ric(h_0) + \frac{1}{\tau}h_0 -
\mathcal{L}_{\psi_*\frac{d}{dt}\psi^{-1}}h_0,
\end{equation}
From $\psi\circ\psi^{-1} = \Id$, by taking a time derivative, we see
that $\psi^*\frac{d}{dt}\psi^{-1} +
\psi^*\mathcal{L}_{\frac{d}{dt}\psi}\psi^{-1} = 0$ and since $\psi$ is 
a diffeomorphism, we get that 
\begin{equation}
\label{equation-equation_psi_inv}
\frac{d}{dt}\psi^{-1} = -\mathcal{L}_{\frac{d}{dt}\psi}\psi^{-1} =
\mathcal{L}_{\nabla f(\psi)}\psi^{-1}.
\end{equation}
Since $\{f(t)\}_{0\le t<\infty}$ are the minimizers for $\mathcal{W}$,
there are uniform $C^{k+2,\alpha}$ estimates on $f(t)$. Since
$\frac{d}{dt}\psi = -\nabla f(\psi)$, there are uniform
$C^{k+1,\alpha}$ bounds on $\psi$, for $t\in [0,B]$.  This together
with (\ref{equation-equation_psi_inv}) yields
$|\psi^{-1}|_{k,\alpha}\le C(B)$, for $t\in [0,B]$. Let $\tilde{g}(t)
= \psi^{-1}g(t)$. Then $\tilde{g}(t)$ satisfies the equation
$$\frac{d}{dt}\tilde{g} = -2\ric(\tilde{g}) + \frac{1}{\tau}\tilde{g}
- \mathcal{L}_{\psi_*\frac{d}{dt}\psi^{-1}}\tilde{g},$$
and
$$|\tilde{g}(t_i+t) - h_0|_{k,\alpha} \le |\psi^{-1}||g(t_i+t) - h(t)|
\le C(B)|g(t_i+t) - h(t)| \to 0,$$ 
when $i\to\infty$, uniformly on
$M\times[0,B]$ (that implies $\tilde{g}(t_i+t)\to h_0$ uniformly on
compact subsets of $M\times [0,\infty)$). The proof of Proposition
\ref{proposition-proposition_extension}, after minor modifications can
be used to get the following result that tell us how to find an
appropriate gauge in the case of convergence toward the solitons
instead of Einstein metrics.

\begin{theorem}
\label{theorem-theorem_extension}
For every $L>0$ and an integer $k$, there exists
$\epsilon_0=\epsilon_0(L,k)$ such that for every $\epsilon <
\epsilon_0$ we can find $i_0=i_0(L,\epsilon,k)$, so that whenever
$i\ge i_0$ there is a gauge $\phi(t)$ on $M\times [t_i,t_i+L]$ such
that $\phi g$ is in a standard form around $h_0$ (see Remark
\ref{remark-remark_standard_form} below), $|\phi\tilde{g} -
h_0|_{k,\alpha} < \epsilon$ and $|\phi - \Id|_{k,\alpha} < \epsilon$.
\end{theorem}

\begin{definition}
\label{definition-definition_integrability}
A limit soliton $h(0)$ is said to be {\bf integrable} if for
every solution ${\bf a}$ of a linearized deformation equation
$$\frac{d}{du}(\ric_{g_u} +
\mathcal{L}_{\psi_*\frac{d}{dt}\psi^{-1}}g_u -
\frac{1}{\tau}(g_u)_{ij})|_{u=0} = 0,$$ 
with $g_0 = h_0$ there exists a path of solitons $h_u$, satisfying
the soliton equation 
\begin{equation}
\label{equation-equation_int_def5}
\ric_{h_u} + \mathcal{L}_{\psi_*\frac{d}{dt}\psi^{-1}}h_u -
\frac{1}{\tau}(h_u)_{ij} = 0,
\end{equation}
with $u\in (-\epsilon,\epsilon)$ and
$h_{0} = h(0)$ such that
$$\frac{d}{du}|_{u=0}h_u = {\bf a}.$$
\end{definition}

\begin{remark}
\label{remark-remark_standard_form}
In the context of Theorem \ref{theorem-theorem_uniqueness}, to say that
$\bar{g}(t)$ is in a {\bf standard form around $h_0$} means that
$\bar{g}$ satisfies the following equation
\begin{equation}
\label{equation-equation_real}
\frac{d}{dt}\bar{g} = -2\ric(\bar{g}) + \frac{1}{\tau}\bar{g} +
P_{h_0}(\bar{g}) - \mathcal{L}_{\psi_*\frac{d}{dt}\psi^{-1}}\bar{g},
\end{equation}
where $P_{h_0}(\tilde{g}) = \nabla_i V_j + \nabla_j V_i$ and $V^k =
\tilde{g}^{pq}(\Gamma_{pq}^k(\tilde{g}) - \Gamma_{pq}^k(h_0))$. We
will write $h_0$ for $h(0)$ in a further discussion.
\end{remark} 

Choose $i_0$, $\phi$ as in Theorem \ref{theorem-theorem_extension}
with $3L$ instead of $L$. Denote by $||\cdot||_{a,b} =
\int_a^b|\cdot|$, where $|\cdot|$ is just the $L^2$ norm. Let $\pi$
denote an orthogonal projection on the subspace $\ker(-\frac{d}{dt} +
\Delta + \frac{1}{\tau} + U)_{M\times [t_{i_0},t_{i_0}+L]}$, with
respect to norm $||\cdot||_{t_{i_0},t_{i_0}+L}$, where $U$ is a linear
first-order expression that comes out after linearizing the equation
(\ref{equation-equation_real}). Let $g_1$ be a suitable chosen
soliton. Denote by $k = \phi g - g_1$ and put $\pi k = (\pi
k)_{\uparrow} + (\pi k)_{\downarrow} + (\pi k)_0$. The integrability
assumption on $h_0$ enters when we choose $g_1$ so that $(\pi k)_0 =
0$. Look at the explanation for $(\cdot)_{\uparrow}$,
$(\cdot)_{\downarrow}$ and $(\cdot)_0$, just after the equation
(\ref{equation-equation_lin1}) below.

\begin{lemma}
\label{lemma-lemma_integrability}
Let $h_0$ be an integrable limit soliton. Then if $\tau < \tau(n,L)$,
for any cylinder $M\times[t_{i_0},t_{i_0}+L]$ there is a soliton $g_1$
satisfying $P_{h_0}(g_1) = 0$ and equation
(\ref{equation-equation_int_def5}), and such that $(\pi k)_0 =
0$. Moreover, if
$$\sup_{[t_{i_0},t_{i_0}+L]}|\phi g(t) - h_0| < \tau,$$
then
\begin{equation}
\label{equation-equation_cond_comp}
||g_1 - h_0||_{t_{i_0},t_{i_0}+L} \le
2||\pi(\phi g(t) - h_0)||_{t_{i_0},t_{i_0}+L}.
\end{equation}
\end{lemma}

\begin{proof}
The proof of this lemma follows the proof of Lemma $5.56$ in
\cite{cheeger1994}. The integrability assumption implies that the set
of metrics $\tilde{g}$ satisfying
$$\ric(\tilde{g}) - \frac{1}{\tau}\tilde{g} +
\mathcal{L}_{\psi_*\frac{d}{dt}\psi^{-1}}\tilde{g} = 0,$$
$$P_{h_0}(\tilde{g}) = 0,$$ 
has a natural smooth manifold structure
near $h_0$. Let $\mathcal{V}$ be a sufficiently small Euclidean
neighborhood of $h_0$. The tangent space to $\mathcal{V}$ at $h_0$
is naturally identified with
$$\mathcal{K} = \{a\in\ker(-\frac{d}{dt} + \Delta + \frac{1}{\tau}
+ U)|P_{h_0}a = 0\}.$$ Define $\mathcal{\psi}:
\mathcal{V}\to\mathcal{K}$ by
$$\mathcal{\psi}(\tilde{g}) = \sum_i\langle\tilde{g},B_i\rangle B_i,$$
where $B_i$ is an orthonormal basis for $\mathcal{K}$ with respect to
a natural inner product. $\mathcal{\psi}$ is a smooth map and the
differential of $\mathcal{\psi}$ is the identity map. We can use now
the implicit function theorem and Lemma
\ref{lemma-lemma_regularity_est} to finish the proof of the Lemma
\ref{lemma-lemma_integrability}.
\end{proof}

The inequality (\ref{equation-equation_cond_comp}) implies that
$|g_1-h_0| \le 2\sup_{[t_{i_0},t_{i_0}+L]}|\pi(\phi g(t) - h_0|$,
where $|\cdot|$ is just the usual $L^2$ norm. The linearization of the
right hand side of the equation $\frac{d}{dt}\phi g = Q(\phi g)$,
satisfied by $\phi g$, where $\phi$ is a gauge chosen as in Theorem
\ref{theorem-theorem_extension} is $DQ(k) = \Delta k + \frac{1}{\tau}k
+ U$, where $U$ is a linear first-order expression in $k$ and a
Laplacian and $U$ are with respect to metric $\phi g$. Let $F$ be a
solution of
\begin{equation}
\label{equation-equation_lin1}
\frac{d}{dt}F = \mathcal{L}F,
\end{equation} 
where $\mathcal{L} = \Delta + \frac{1}{\tau} + U$ and the Laplacian
and $U$ are this time given with a respect to a fixed metric (in our
case we will take metric $h_0$). Let $\{\lambda_k\}$ be the set of
eigenvalues of $\mathcal{L}$. We can write $F = F_{\uparrow} +
F_{\downarrow} + F_0$, where $F_{\uparrow}(t) =
\sum_{\lambda_k<0}a_ke^{-\lambda_kt}$, $F_{\downarrow}(t) =
\sum_{\lambda_k>0}a_k e^{-\lambda_kt}$, and $F_0$ is a projection of
$F$ to a kernel of $\mathcal{L}$.

The basic parabolic estimates (for example similarly as in
\cite{simon1983} and \cite{cheeger1994}) yield the following.

\begin{lemma}
\label{lemma-lemma_regularity_est}
There exists $\tau > 0$ such that for any solution $\eta$ of
(\ref{equation-equation_lin1}) with $|g_1-h_0|_{k+2,\alpha} \le \tau$,
we have that
$$\sup_{(t_0,t_0+L)}|\eta|_{k,\alpha} \le C\sup_{(t_0,t_0+L)}|\eta|,$$
where the first norm is $C^{k,\alpha}$ norm and the last norm is $L^2$
norm.
\end{lemma}

\begin{lemma}
\label{lemma-lemma_alpha}
There exists $\alpha > 1$ such that
\begin{equation}
\label{equation-equation_growth}
\sup_{[L,2L]}|F_{\uparrow}| \ge \alpha\sup_{[0,L]}|F_{\uparrow}|,
\end{equation}
\begin{equation}
\label{equation-equation_decay}
\sup_{[L,2L]}|F_{\downarrow}| \le \alpha^{-1}\sup_{[0,L]}|F_{\downarrow}|.
\end{equation}
The norms considered above are standard $L^2$ norms.
\end{lemma}

\begin{proof}
We will prove only (\ref{equation-equation_growth}), since the proof
of (\ref{equation-equation_decay}) is similar. Let $\delta =
\min\{|\lambda_k| \neq 0\} > 0$.
\begin{eqnarray*}
\sup_{[L,2L]}|F_{\uparrow}| - \alpha\sup_{[0,L]}|F_{\uparrow}| &=&
\sup_{[0,L]}\sum_{\lambda_k<0}a_k^2e^{-2\lambda_kt}e^{-2\lambda_kL}
- \alpha\sup_{[0,L]}\sum_{\lambda_k<0}a_k^2e^{-2\lambda_kt} \\ 
&\ge& \sup_{[0,L]}\sum_{\lambda_k<0}a_k^2e^{-2\lambda_kt}(e^{2\delta L} -
\alpha),
\end{eqnarray*}
which is positive, if $e^{2\delta L} > \alpha$. We can choose 
$\alpha = e^{\delta L} > 1$.
\end{proof}

\begin{lemma}
\label{lemma-lemma_beta}
There exists $\beta < \alpha$ such that if
\begin{equation}
\label{equation-equation_impl1}
\sup_{[L,2L]}|F| \ge \beta\sup_{[0,L]}|F|,
\end{equation}
then
\begin{equation}
\label{equation-equation_impl2}
\sup_{[2L,3L]}|F| \ge \beta\sup_{[L,2L]}|F|,
\end{equation}
and if
\begin{equation}
\label{equation-equation_impl3}
\sup_{[2L,3L]}|F| \le \beta^{-1}\sup_{[L,2L]}|F|,
\end{equation}
then
\begin{equation}
\label{equation-equation_impl4}
\sup_{[L,2L]}|F| \le \beta^{-1}\sup_{[0,L]}|F|.
\end{equation}
Moreover, if $F_0 = 0$ at least one of (\ref{equation-equation_impl2}),
(\ref{equation-equation_impl4}) holds.
\end{lemma}
The proof of Lemma \ref{lemma-lemma_beta} is almost the same to the
proof of analogous lemma ($5.31$) in \cite{cheeger1994}. We can choose
$\beta$ to be of order $e^{\frac{L\delta}{4}}$.

Let $\eta = \phi g - g_1$, where $\phi$ is chosen as in Theorem
\ref{theorem-theorem_extension} and $g_1$ is a soliton as in Lemma
\ref{lemma-lemma_integrability} which does not depend on $t$ for a 
considered time interval of length $L$.

\begin{lemma}
\label{lemma-lemma_linear_equation}
\begin{equation}
\label{equation-equation_linear_eq}
\frac{d}{dt}(\phi g - g_1) = \Delta_{h_0}(\phi g - g_1) +
\frac{1}{\tau}(\phi g - g_1) + F(\phi g,h_0,g_1) + U(\phi g - g_1)
\end{equation}
where $|F(\phi g,h,g_1)|_{k,\alpha} \le C(|g_1-h_0| +
|\eta|_{k,\alpha})|\nabla^2\eta|_{k-2,\alpha} +
C(|\nabla(g_1-h_0)|_{k-1,\alpha} + |\nabla \eta|_{k-1,\alpha})
|\nabla\eta|_{k-1,\alpha}$ and $U$
is a first order linear expression in $\phi g - g_1$.
\end{lemma}

\begin{proof}
Since both $\phi g$ and $g_1$ are in a standard form around $h_0$
(recall that $P_{h_0}(g_1) = 0$), by using a formula for linearization
of a second order operator $-2\ric(\phi g) + P_{h_0}(\phi g)$, we get
\begin{eqnarray}
\label{equation-equation_est_equ}
\frac{d}{dt}(\phi g - g_1) &=& (-2\ric(\phi g) + P_{h_0}(\phi g) -
\mathcal{L}_{\phi_*\psi_*\frac{d}{dt}\psi^{-1}}\phi g) - \\ 
&-& (-2\ric(g_1) + P_{h_0}(g_1) - 
\mathcal{L}_{\psi_*\frac{d}{dt}\psi}g_1) +  
\frac{1}{\tau}(\phi g - g_1) \nonumber \\
&=& \Delta_{\phi g}(\phi g - g_1) + \frac{1}{\tau}(\phi g - g_1) + 
U(\phi g - g_1) + \tilde{F}(\phi g,g_1), \nonumber
\end{eqnarray}
where $|\tilde{F}(\phi g,g_1)|_{k,\alpha} \le
C(|\eta|_{k,\alpha}|\nabla^2\eta|_{k-2,\alpha} +
|\nabla\eta|^2_{k,\alpha})$, by a similar computation to a computation
in \cite{cheeger1994}. Furthermore, $\Delta_{\phi g}\eta =
\Delta_{h_0}\eta + (\Delta_{\phi g} - \Delta_{h_0})\eta$ and since
$|\phi g - h_0|_{k,\alpha} \le C(|\eta|_{k,\alpha} + |g_1 -
h_0|_{k,\alpha})$, we have that $|(\Delta_{\phi g} - \Delta_{h_0})\eta|
\le C(|\eta|_{k,\alpha} + |g_1 -
h_0|_{k,\alpha})|\nabla^2\eta|_{k,\alpha}$. The Lemma
\ref{lemma-lemma_linear_equation} now follows.
\end{proof}

We assume that $|g_1-h_0|_{k,\alpha} < \epsilon$. Let $k$ be a
solution to (\ref{equation-equation_linear_eq}). Then we have the
following Proposition.
 
\begin{proposition}
\label{proposition-proposition_comparison}
There exists $\epsilon_0 > 0$, depending on the uniform bounds on the
geometries $g(t)$, such that if $\epsilon < \epsilon_0$, then if
\begin{equation}
\label{equation-equation_sol1}
\sup_{[L,2L]}|k| \ge \beta\sup_{[0,L]}|k|,
\end{equation}
then  
\begin{equation}
\label{equation-equation_sol2}
\sup_{[2L,3L]}|k| \ge \beta\sup_{[L,2L]}|k|,
\end{equation}
and if
\begin{equation}
\label{equation-equation_sol3}
\sup_{[2L,3L]}|k| \le \beta^{-1}\sup_{[L,2L]}|k|,
\end{equation}
then
\begin{equation}
\label{equation-equation_sol4}
\sup_{[L,2L]}|k| \le \beta^{-1}\sup_{[0,L]}|k|,
\end{equation}
-Moreover, if $(\pi k)_0 = 0$, at least one of
(\ref{equation-equation_sol2}), (\ref{equation-equation_sol4}) holds.
\end{proposition}

\begin{proof}
Assume there exist a sequence of gauges $\phi_i$ and constants
$\tau_i\to 0$, such that $|\eta_i|_{k,\alpha} = |\phi_i g -
h|_{k,\alpha} \le \tau_i \to 0$, but for which none of the assertions
in Proposition \ref{proposition-proposition_comparison} holds. Let
$\psi_i = \frac{\eta_i}{\sup_{[L,2L]}|\eta_i|}$. Then in view of Lemma
\ref{lemma-lemma_regularity_est}, from standard compactness results
(as in \cite{cheeger1994}) we get that for a subsequence $\psi_i
\stackrel{C^{k,\alpha}}{\to} \psi$ and

$$\frac{d}{dt}\psi = \Delta_h \psi + U(\psi) + \frac{1}{\tau}\psi,$$
where $\psi$ has a property that contradicts Lemma
\ref{lemma-lemma_beta}. Recall that $\beta$ is of order
$e^{\frac{\epsilon L}{4}}$.
\end{proof}

\begin{proof}[Proof of Theorem \ref{theorem-theorem_uniqueness}]
We will adopt the notation from above. Take some $L>0$ big enough (we
will see later how big we want to make it) and choose $\epsilon_0 > 0$
as in Theorem \ref{theorem-theorem_extension} so that the Theorem
holds for $\epsilon_0$, and $3L$. For every $\epsilon < \epsilon_0$
there exists $i_0$ such that for every $i\ge i_0$ there exists a gauge
$\phi$ so that $\phi$ satisfies all the conditions in Theorem
\ref{theorem-theorem_extension}, that is $\phi g$ is in a standard
form around $h_0$, $|\phi g - h_0|_{k,\alpha} < \epsilon$ and
$|\frac{d}{dt}\phi g| < \tilde{\epsilon}$ on $M\times [t_i,t_i+3L]$,
where $\tilde{\epsilon}$ is comparable to $\epsilon$. For each $t_{i}$
pick up the largest possible $L'$ (we will omit emphasizing a
dependence of $L'$ on $i$ and we will call it just $L'$, since it is
irrelevant for further discussion) such that {\bf (**)} $\phi$ is
defined on $M\times [t_i,t_i+L')$, $\phi g$ is in a standard form
around $h$ and $|\phi g - h|_{k,\alpha} < \epsilon$ and
$\min_{[t_i,t_i+3L]}|\phi g - h|_{k,\alpha} <
\frac{\epsilon}{1000}$. Divide $[t_i,t_i+L')$ into the subintervals of
length $L$ and assume that $N$ is the largest number such that
$[t_i+(N-1)L,t_i+NL]\subset [t_i,t_i+L')$.

Notice that for $L$ chosen above, from the proof of Theorem
\ref{theorem-theorem_extension}, all the estimates that we have got on
$|\phi - \Id|_{k,\alpha}$ in the previous subsection depend on a
polynomial in $L$ (call it $q(L)$), whose coefficients depend only on
a dimension, an integer $k$ and the uniform bounds on geometries
$g(t)$. By the estimates established in Proposition
\ref{proposition-proposition_extension}, we can increase $i_0$ if
necessary, so that
\begin{enumerate}
\item
\label{item-item_0}
For every $i\ge i_0$ we can find a gauge on $M\times[t_i,t_i+3L]$,
such that $\sup_{[t_i,t_i+3L]}|\phi g(t) - h_0|_{k,\alpha} <
\frac{\epsilon}{1000e^{L\delta/4}}$.
\item
\label{item-item_1}
If the initial data $\phi(s)$ is such that $|\phi(s) - \Id|_{k,\alpha}
< \frac{\epsilon}{e^{L\delta/8}}$ and $|\phi(s)^*g(s)-h_0|_{k,\alpha}
< \frac{\epsilon}{e^{L\delta/8}}$, where $s\in [t_i,t_i+L']$, for
$i\ge i_0$, then $\phi$ can be extended to interval $[s,s+3L]$ such
that $\sup_{[s,s+3L]}|\phi g - h_0|_{k,\alpha} <
\frac{\epsilon}{100p(L)}$ (we might need increase $i_0$ for this to
hold). Polynomial $p(L)$ can be any polynomial with leading
coefficient $1$ and with a degree that is e.g. one more than a degree
of $q(L)$.
\item
\label{item-item_2}
If the initial data is such that $|\phi(s) - \Id|_{k,\alpha} <
\frac{\epsilon}{p(L)}$ and $|\phi(s)^*g(s)-h_0|_{k,\alpha} <
\frac{\epsilon}{p(L)}$, where $s\in [t_i,t_i+L']$, for $i\ge i_0$,
then $\phi$ can be extended on interval $[s,s+3L]$ such that
$\sup_{[s,s+3L]}|\phi g - h_0|_{k,\alpha} < \epsilon$.
\end{enumerate}

We want to show that there exists $i$ (for sufficiently big $L$, so
that above holds) such that a corresponding $L' = \infty$. Assume that
for all $i\ge i_0$ and all $\epsilon > 0$, $L'<\infty$. Denote by $I_j
= [t_i+jL,t_i+jL+L]$. Assume that $\epsilon$ is small enough so that
we can apply Lemma \ref{lemma-lemma_integrability}, that is for every
$j$ there exists a soliton $g_j$ such that $(\pi(\phi g - g_j))_0 = 0$
on $I_j$ and therefore by Proposition
\ref{proposition-proposition_comparison}, $\phi g - g_j$ either
satisfies a growth condition ((\ref{equation-equation_sol1})
$\Rightarrow$ (\ref{equation-equation_sol2})) or a decay condition
((\ref{equation-equation_sol3}) $\Rightarrow$
(\ref{equation-equation_sol4})). Moreover, $|g_j-h_0| \le
2\sup_{I_j}|\pi(\phi g - h_0)| \le C\sup_{I_j}|\phi g -
h_0|_{k,\alpha}$. We need to consider two cases.

\begin{case}
Assume that for all $i_0$ and all $i\ge i_0$, where $i_0 = i_0(L)$ is
chosen as in Theorem \ref{theorem-theorem_extension} for $L$ big
enough (so that (\ref{item-item_0}), (\ref{item-item_1}) and
(\ref{item-item_2}) hold), and for all the intervals $I_j$ (that are
defined with respect to $t_i$; we want to omit double indices) for
which we have $\sup_{I_j}|\phi g - h_0|_{k,\alpha} \le
\frac{\epsilon}{100p(L)}$, $\phi g - g_j$ satisfies a decay condition
on $I_j$ (recall that $L^2$ norms are considered in a growth and a
decay condition).
\end{case}

By using Proposition \ref{proposition-proposition_comparison}
inductively, we get that 
$$\sup_{I_l}|\phi g - g_j| \le \frac{1}{\beta^l}\sup_{I_1}|\phi g -
g_j|,$$ 
for all $l\le j$. Moreover, $\sup_{I_1}|\phi g - g_j| \le
\sup_{I_1}|\phi g - h_0| + |g_j - h_0| \le \sup_{I_1}|\phi g - h_0| +
2\sup_{I_j}|\phi g - h_0| < \frac{3\epsilon}{100p(L)}$, which yields 
$$\sup_{I_l}|\phi g - g_j| \le
\frac{1}{\beta^l}\frac{3\epsilon}{100p(L)},$$ 
By Lemma \ref{lemma-lemma_regularity_est} we may assume that $\sup_{I_l}|\phi g
- g_j|_{k+2,\alpha} \le \frac{\epsilon}{\beta^l}$. Whenever we increase
$L$ (the necessity for $L$ being increased will depend only on the
uniform estimates), we can choose an appropriate $\epsilon_0$ as in
Theorem \ref{theorem-theorem_extension} and take any $\epsilon <
\epsilon_0$. Each time we do that we might have to increase $i_0$
(depending on $\epsilon < \epsilon_0$). Therefore, on $M\times
I_l$, for $l\le j$ we have

\begin{eqnarray*}
|\frac{d}{dt}\phi g|_{k,\alpha} &=& |\frac{d}{dt}(\phi g - g_j)|_{k,\alpha} \\
&=& (-2\ric(\phi g) + 2\ric(g_j)) + \frac{1}{\tau}(\phi g - g_j)
+ (P_{h_0}(\phi g) - P_{h_0}(g_j)) + \mathcal{L}_{\psi_*\frac{d}{dt}\psi^{-1}}
(g_j - \phi g) \\
&\le& C\sup_{I_l}|\phi g - g_j|_{k+2,\alpha} < \frac{C\epsilon}{\beta^l}.
\end{eqnarray*}
For every $l \le j$, since $\frac{d}{dt}\phi g = \frac{d}{dt}(\phi g - h_0)$,
we have that
\begin{eqnarray*}
\sup_{I_l}|\phi g - h_0|_{k,\alpha} &\le& 2L\sup_{I_l\cup I_{l-1}}
|\frac{d}{dt}\phi g|_{k,\alpha} + \sup_{I_{l-1}}|\phi g -
h_0|_{k,\alpha} \\ 
&\le& 2LC\frac{\epsilon}{\beta^{l-1}} +
2LC\frac{\epsilon}{\beta^{l-2}} + \dots + 2LC\frac{\epsilon}{\beta} +
\sup_{I_2}|\phi g - h_0|_{k,\alpha} \\ 
&\le& \sup_{I_2}|\phi g - h_0|_{k,\alpha} + \frac{2LC\epsilon}{\beta - 1},
\end{eqnarray*}
which can be made smaller than $\frac{\epsilon}{e^{L\delta/8}}$ for
$L$ chosen big enough at the beginning. By condition
\ref{item-item_1}, for big values of $i$ we can extend $\phi$ on
$I_{j+1}$ so that $\sup_{I_{j+1}}|\phi g - h_0| <
\frac{\epsilon}{100p(L)}$ and it has to coincide with our previously
constructed $\phi$ on $I_{j+1}$. We can continue a described procedure
by looking now at intervals $I_{j}$ and $I_{j+1}$ replaced by
intervals $I_{j+1}$ and $I_{j+2}$ respectively. If we repeat this
sufficiently many times, we will reach the interval $I_{N-1}$ with
$$\sup_{I_{N-1}}|\phi g - h_0|_{k,\alpha} <
\frac{\epsilon}{100p(L)}.$$ 
By condition \ref{item-item_2} we will now
be able to extend $\phi$ (for sufficiently big values of $i$) to
interval $[t_i+(N-1),t_i+(N+1)L]$, with
$\sup_{[t_i+(N-1),t_i+(N+1)L]}|\phi g - h_0| < \epsilon$ holding.
Since $(N+1)L > L'$, this estimate contradicts a maximality of $L'$
with properties {\bf (**)}. Therefore, either there exists $i$ such
that a corresponding $L' = \infty$, or we have a following case
holding.

\begin{case}
There are some $L$, $i$ and $j$ for which $\sup_{I_j}|\phi g -
h_0|_{k,\alpha} < \frac{\epsilon}{100p(L)}$, and $\phi g - h_0$
satisfies a growth condition on $I_j$ ($I_j$ is defined with respect
to $t_i$).
\end{case}

By using Proposition \ref{proposition-proposition_comparison}
inductively, we would have that

\begin{eqnarray*}
\sup_{I_{N-1}}|\phi g - g_j| &<& \frac{1}{\beta}\sup_{I_N}|\phi g - g_j| \\
&\le& \frac{1}{\beta}(\sup_{I_N}|\phi g - h_0| + |g_j - h_0|) \\
&\le& \frac{1}{\beta}(\sup_{I_N}|\phi g - h_0| + 2\sup_{I_j}|\phi g - h_0|) \\
&<& \frac{3\epsilon}{\beta}.
\end{eqnarray*}
Moreover, if we use Lemma \ref{lemma-lemma_regularity_est}, together with the
estimate
\begin{eqnarray*}
\sup_{I_{N-1}}|\phi g - h_0| &\le& \sup_{I_{N-1}}|\phi g -
g_j| + |g_j - h_0| \\
&\le& \sup_{I_{N-1}}|\phi g - g_j| + 2\sup_{I_j}|\phi g - h_0| \\
&<& \frac{3\epsilon}{\beta} + \frac{\epsilon}{100p(L)},
\end{eqnarray*}
which can be made smaller than $\frac{\epsilon}{p(L)}$, by condition
\ref{item-item_2} we can extend $\phi$ to an interval
$[t_i+(N-1)L,t_i+(N+1)L]$ (if $i$ is big enough), with $|\phi g -
h_0|_{k,\alpha} < \epsilon$ holding. We again get a contradiction as
in the previous case if we assume $L' < \infty$ for all $i$.

Therefore, there exists $i_0$ such that a gauge $\phi$ can be
constructed on $M\times [t_{i_0},t_{i_0}+L')$, satisfying properties
{\bf (**)} and such that a corresponding $L' = \infty$. Consider again
$I_j = [t_{i_0}+jL,t_{i_0}+jL+L]$ and the corresponding $g_j$ that are
found by Lemma \ref{lemma-lemma_integrability}, such that for $k_j =
\phi g - g_j$ we have that $(\pi k_j)_0 = 0$ on $M\times I_j$. Notice
that a decay condition ((\ref{equation-equation_sol3} $\Rightarrow$
(\ref{equation-equation_sol4})) holds for all $j$. If there existed
some $j$ for which it were not true, by using Proposition
\ref{proposition-proposition_comparison} inductively and standard
parabolic estimates (Lemma \ref{lemma-lemma_regularity_est}), we would
find that

$$\epsilon > \sup_{[t_{i_0}+(N-1)L,t_{i_0}+NL]}|\phi g - g_j| \ge
\beta^{N-j}\sup_{I_j}|\phi g - g_j|,$$ 
for all $N$ and we
would get a contradiction by letting $N$ tend to infinity (if
$\sup_{I_j}|\phi g - g_j| = 0$, our metric $g(t_{i_0}+jL)$
would be a soliton satisfying (\ref{equation-equation_soliton_h0}) and
it would stay so for all later times which is not an interesting
case). This means we have a decay for all times if we do not start with
a soliton. 

After passing to a subsequence, we may assume that for some metric
$g_{\infty}$ that satisfies a soliton type equation
$\lim_{j\to\infty}|g_{j_p}-g_{\infty}|_{k,\alpha'} = 0$, where $\alpha' < 
\alpha$.

\begin{claim}
\label{claim-claim_decay}
$\lim_{p\to\infty}\sup_{I_{j_p}}|k_{j_p}|_{k,\alpha} = 0$.
\end{claim}

\begin{proof}
If it were not the case, there would exist a subsequence of $j_p$
(denote it by the same symbol) such that $\phi g - g_{j_p}$ would
satisfy a growth condition, that is
$$\sup_{[t_{i_0}+(N-1)L,t_{i_0}+NL]}|\phi g - g_{j_p}| \ge
\beta^{(N-j_p)}\sup_{I_{j_p}}|\phi g - g_{j_p}|,$$ 
for all $N$, where $\beta$ can be taken to be $e^{\frac{L\delta}{4}}$ 
and by taking $N\to\infty$ we immediately get a contradiction, since
$\sup_{[t_{i_0}+(N-1)L,t_{i_0}+NL]}|\phi g - g_{j_p}| < C\epsilon$.
\end{proof}

As in the proof of the claim above, we get that $\phi g - g_{j_p}$ has
to satisfy a decay condition for all $p$. By Claim
\ref{claim-claim_decay}, by using Proposition
\ref{proposition-proposition_comparison} inductively and by standard
parabolic estimates (Lemma \ref{lemma-lemma_regularity_est}) we find
that for some $c>0$,
$$|\phi g - g_{\infty}|_{k,\alpha} \le
ce^{-\frac{\delta(t-t_{i_0})}{4}},$$ 
for $t\in [t_{i_0}+(N-1)L,t_{i_0}+NL]$ and for all $N > 0$, that is
\begin{equation}
\label{equation-equation_exp_decay}
|\phi g - g_{\infty}|_{k,\alpha} \le Ce^{-ct},
\end{equation}
for all $t\ge t_{i_0}$. (\ref{equation-equation_exp_decay}) implies
that $|g(t) - \phi^{-1}g_{\infty}|_{C^0} <
Ce^{-ct}$. $\phi^{-1}g_{\infty}$ is a soliton that moves by
diffeomorphisms $\phi(t)^{-1}$ and therefore is determined by metric
$\phi^{-1}(t_{i_0})g_{\infty}$. Since $h_0$ is a limit soliton of
metrics $g(t_i)$, $h_0$ and $\phi^{-1}(t_{i_0})g_{\infty}$ differ only
by a diffeomorphism, that is $\eta\phi^{-1}(t_{i_0})g_{\infty} = h_0$
for some diffeomorphism $\eta$. Let finally $\phi' =
\eta\phi^{-1}(t_{i_0})\phi(t)$. Then,
$$|\phi'g(t) - h_0|_{k,\alpha} < Ce^{-ct},$$
that is $\phi'g(t)$ converges to a soliton $h_0$ exponentially as
$t\to\infty$. We know that $h(t) = \psi(t) h_0$ and therefore,
$$|\psi\phi'g(t) - h(t)|_{C^0} \le Ce^{-ct}.$$ 
This finishes the proof of Theorem \ref{theorem-theorem_uniqueness}.
\end{proof}

\end{subsection}

\end{section}

\end{document}